%
\ifx\shlhetal\undefinedcontrolsequence\let\shlhetal\relax\fi
\input amstex
\NoBlackBoxes
\documentstyle {amsppt}
\expandafter\ifx\csname bib4plain.tex\endcsname\relax
  \expandafter\gdef\csname bib4plain.tex\endcsname{}
\else \message{Hey!  Apparently you were trying to \string twice.   This does not make sense.}
\errmessage{Please edit your file (probably \jobname.tex) and remove
any duplicate ``\string\input'' lines} \fi

\def\renewcommand{\newcommand}	       
\edef\cite{\the\catcode`@}%
\catcode`@ = 11
\let\@oldatcatcode = \cite
\chardef\@letter = 11
\chardef\@other = 12
%
%
%
%
\def\@innerdef#1#2{\edef#1{\expandafter\noexpand\csname #2\endcsname}}%
%
%
\@innerdef\@innernewcount{newcount}%
\@innerdef\@innernewdimen{newdimen}%
\@innerdef\@innernewif{newif}%
\@innerdef\@innernewwrite{newwrite}%
%
%
%
\def\@gobble#1{}%
%
%
%
\ifx\inputlineno\@undefined
   \let\@linenumber = \empty 
\else
   \def\@linenumber{\the\inputlineno:\space}%
\fi
%
%
%
\def\@futurenonspacelet#1{\def\cs{#1}%
   \afterassignment\@stepone\let\@nexttoken=
}%
\begingroup 
\def\\{\global\let\@stoken= }%
\\ 
\endgroup
\def\@stepone{\expandafter\futurelet\cs\@steptwo}%
\def\@steptwo{\expandafter\ifx\cs\@stoken\let\@@next=\@stepthree
   \else\let\@@next=\@nexttoken\fi \@@next}%
\def\@stepthree{\afterassignment\@stepone\let\@@next= }%
%
%
%
\def\@getoptionalarg#1{%
   \let\@optionaltemp = #1%
   \let\@optionalnext = \relax
   \@futurenonspacelet\@optionalnext\@bracketcheck
}%
%
%
\def\@bracketcheck{%
   \ifx [\@optionalnext
      \expandafter\@@getoptionalarg
   \else
      \let\@optionalarg = \empty
      \expandafter\@optionaltemp
   \fi
}%
\def\@@getoptionalarg[#1]{%
   \def\@optionalarg{#1}%
   \@optionaltemp
}%
%
%
%
\def\@nnil{\@nil}%
\def\@fornoop#1\@@#2#3{}%
\def\@for#1:=#2\do#3{%
   \edef\@fortmp{#2}%
   \ifx\@fortmp\empty \else
      \expandafter\@forloop#2,\@nil,\@nil\@@#1{#3}%
   \fi
}%
\def\@forloop#1,#2,#3\@@#4#5{\def#4{#1}\ifx #4\@nnil \else
       #5\def#4{#2}\ifx #4\@nnil \else#5\@iforloop #3\@@#4{#5}\fi\fi
}%
\def\@iforloop#1,#2\@@#3#4{\def#3{#1}\ifx #3\@nnil
       \let\@nextwhile=\@fornoop \else
      #4\relax\let\@nextwhile=\@iforloop\fi\@nextwhile#2\@@#3{#4}%
}%
%
%
%
\@innernewif\if@fileexists
\def\@testfileexistence{\@getoptionalarg\@finishtestfileexistence}%
\def\@finishtestfileexistence#1{%
   \begingroup
      \def\extension{#1}%
      \immediate\openin0 =
         \ifx\@optionalarg\empty\jobname\else\@optionalarg\fi
         \ifx\extension\empty \else .#1\fi
         \space
      \ifeof 0
         \global\@fileexistsfalse
      \else
         \global\@fileexiststrue
      \fi
      \immediate\closein0
   \endgroup
}%
%
%
%
%
\def\bibliographystyle#1{%
   \@readauxfile
   \@writeaux{\string\bibstyle{#1}}%
}%
\let\bibstyle = \@gobble
%
%
\let\bblfilebasename = \jobname
\def\bibliography#1{%
   \@readauxfile
   \@writeaux{\string\bibdata{#1}}%
   \@testfileexistence[\bblfilebasename]{bbl}%
   \if@fileexists
      \nobreak
      \@readbblfile
   \fi
}%
\let\bibdata = \@gobble
%
%
\def\nocite#1{%
   \@readauxfile
   \@writeaux{\string\citation{#1}}%
}%
\@innernewif\if@notfirstcitation
%
%
\def\cite{\@getoptionalarg\@cite}%
%
%
\def\@cite#1{%
   \let\@citenotetext = \@optionalarg
   \printcitestart
   \nocite{#1}%
   \@notfirstcitationfalse
   \@for \@citation :=#1\do
   {%
      \expandafter\@onecitation\@citation\@@
   }%
   \ifx\empty\@citenotetext\else
      \printcitenote{\@citenotetext}%
   \fi
   \printcitefinish
}%
\def\@onecitation#1\@@{%
   \if@notfirstcitation
      \printbetweencitations
   \fi
   \expandafter \ifx \csname\@citelabel{#1}\endcsname \relax
      \if@citewarning
         \message{\@linenumber Undefined citation `#1'.}%
      \fi
      \expandafter\gdef\csname\@citelabel{#1}\endcsname{%
\strut
\vadjust{\vskip-\dp\strutbox
\vbox to 0pt{\vss\parindent0cm \leftskip=\hsize 
\advance\leftskip3mm
\advance\hsize 4cm\strut\openup-4pt 
\rightskip 0cm plus 1cm minus 0.5cm ?  #1 ?\strut}}
         {\tt
            \escapechar = -1
            \nobreak\hskip0pt
            \expandafter\string\csname#1\endcsname
            \nobreak\hskip0pt
         }%
      }%
   \fi
   \csname\@citelabel{#1}\endcsname
   \@notfirstcitationtrue
}%
%
%
\def\@citelabel#1{b@#1}%
%
%
\def\@citedef#1#2{\expandafter\gdef\csname\@citelabel{#1}\endcsname{#2}}%
%
%
%
\def\@readbblfile{%
   \ifx\@itemnum\@undefined
      \@innernewcount\@itemnum
   \fi
   \begingroup
      \def\begin##1##2{%
         \setbox0 = \hbox{\biblabelcontents{##2}}%
         \biblabelwidth = \wd0
      }%
      \def\end##1{}
      %
      %
      \@itemnum = 0
      \def\bibitem{\@getoptionalarg\@bibitem}%
      \def\@bibitem{%
         \ifx\@optionalarg\empty
            \expandafter\@numberedbibitem
         \else
            \expandafter\@alphabibitem
         \fi
      }%
      \def\@alphabibitem##1{%
         \expandafter \xdef\csname\@citelabel{##1}\endcsname {\@optionalarg}%
         \ifx\biblabelprecontents\@undefined
            \let\biblabelprecontents = \relax
         \fi
         \ifx\biblabelpostcontents\@undefined
            \let\biblabelpostcontents = \hss
         \fi
         \@finishbibitem{##1}%
      }%
      \def\@numberedbibitem##1{%
         \advance\@itemnum by 1
         \expandafter \xdef\csname\@citelabel{##1}\endcsname{\number\@itemnum}%
         \ifx\biblabelprecontents\@undefined
            \let\biblabelprecontents = \hss
         \fi
         \ifx\biblabelpostcontents\@undefined
            \let\biblabelpostcontents = \relax
         \fi
         \@finishbibitem{##1}%
      }%
      \def\@finishbibitem##1{%
         \biblabelprint{\csname\@citelabel{##1}\endcsname}%
         \@writeaux{\string\@citedef{##1}{\csname\@citelabel{##1}\endcsname}}%
         \ignorespaces
      }%
      %
      %
      \let\em = \bblem
      \let\newblock = \bblnewblock
      \let\sc = \bblsc
      \frenchspacing
      \clubpenalty = 4000 \widowpenalty = 4000
      \tolerance = 10000 \hfuzz = .5pt
      \everypar = {\hangindent = \biblabelwidth
                      \advance\hangindent by \biblabelextraspace}%
      \bblrm
      \parskip = 1.5ex plus .5ex minus .5ex
      \biblabelextraspace = .5em
      \bblhook
      \input \bblfilebasename.bbl
   \endgroup
}%
%
%
\@innernewdimen\biblabelwidth
\@innernewdimen\biblabelextraspace
%
%
%
\def\biblabelprint#1{%
   \noindent
   \hbox to \biblabelwidth{%
      \biblabelprecontents
      \biblabelcontents{#1}%
      \biblabelpostcontents
   }%
   \kern\biblabelextraspace
}%
%
%
%
\def\biblabelcontents#1{{\bblrm [#1]}}%
%
%
\def\bblrm{\rm}%
%
%
\def\bblem{\it}%
%
%
\def\bblsc{\ifx\@scfont\@undefined
              \font\@scfont = cmcsc10
           \fi
           \@scfont
}%
%
%
\def\bblnewblock{\hskip .11em plus .33em minus .07em }%
%
%
\let\bblhook = \empty
%
%
%
\def\printcitestart{[}
\def\printcitefinish{]}
\def\printbetweencitations{, }
\def\printcitenote#1{, #1}
%
%
%
\let\citation = \@gobble
%
%
%
\@innernewcount\@numparams
%
%
\def\newcommand#1{%
   \def\@commandname{#1}%
   \@getoptionalarg\@continuenewcommand
}%
%
%
\def\@continuenewcommand{%
   \@numparams = \ifx\@optionalarg\empty 0\else\@optionalarg \fi \relax
   \@newcommand
}%
%
%
\def\@newcommand#1{%
   \def\@startdef{\expandafter\edef\@commandname}%
   \ifnum\@numparams=0
      \let\@paramdef = \empty
   \else
      \ifnum\@numparams>9
         \errmessage{\the\@numparams\space is too many parameters}%
      \else
         \ifnum\@numparams<0
            \errmessage{\the\@numparams\space is too few parameters}%
         \else
            \edef\@paramdef{%
               \ifcase\@numparams
                  \empty  No arguments.
               \or ####1%
               \or ####1####2%
               \or ####1####2####3%
               \or ####1####2####3####4%
               \or ####1####2####3####4####5%
               \or ####1####2####3####4####5####6%
               \or ####1####2####3####4####5####6####7%
               \or ####1####2####3####4####5####6####7####8%
               \or ####1####2####3####4####5####6####7####8####9%
               \fi
            }%
         \fi
      \fi
   \fi
   \expandafter\@startdef\@paramdef{#1}%
}%
%
%
%
%
\def\@readauxfile{%
   \if@auxfiledone \else 
      \global\@auxfiledonetrue
      \@testfileexistence{aux}%
      \if@fileexists
         \begingroup
            \endlinechar = -1
            \catcode`@ = 11
            \input \jobname.aux
         \endgroup
      \else
         \message{\@undefinedmessage}%
         \global\@citewarningfalse
      \fi
      \immediate\openout\@auxfile = \jobname.aux
   \fi
}%
%
%
\newif\if@auxfiledone
\ifx\noauxfile\@undefined \else \@auxfiledonetrue\fi
%
%
%
%
\@innernewwrite\@auxfile
\def\@writeaux#1{\ifx\noauxfile\@undefined \write\@auxfile{#1}\fi}%
%
%
%
\ifx\@undefinedmessage\@undefined
   \def\@undefinedmessage{No .aux file; I won't give you warnings about
                          undefined citations.}%
\fi
%
%
\@innernewif\if@citewarning
\ifx\noauxfile\@undefined \@citewarningtrue\fi
%
%
%
\catcode`@ = \@oldatcatcode

\topmatter
\title {{Universal in $(< \lambda)$-Stable Abelian Group} \\
Sh456} \endtitle
\author {Saharon Shelah \thanks {\null\newline
I thank Alice Leonhardt for the beautiful typing \null\newline
Research supported by the United States-Israel Binational Science Foundation;
Pub. No.456} \endthanks} \endauthor
\affil {Institute of Mathematics \\
The Hebrew University \\
Jerusalem, Israel
\medskip
Rutgers University \\
Department of Mathematics \\
New Brunswick, NJ USA} \endaffil

\abstract  A characteristic result is that if $2^{\aleph_0} <\mu < \mu^+ < 
\lambda = \text{ cf}(\lambda) < \mu^{\aleph_0}$, then among the separable
reduced $p$-groups of cardinality $\lambda$ which are
$(<\lambda)$-stable there is no universal one.
\endabstract

\endtopmatter
\document
\newpage

We deal with the existence of reduced separable (abelian) groups of
cardinality $\lambda$ under usual embeddings (or similar classes of abelian
groups or modules).   So this continues Kojman Shelah \cite{KjSh:409},
\cite{KjSh:449}, \cite{KjSh:455}, but we make the presentation self 
contained so we repeat some things from there.  They deal with proving the
non-existence of universal members in cardinality $\lambda$ for various
classes:
\newline
\cite{KjSh:409} deals with linear orders and f.o. theories with the
strict order property; 
\newline
\cite{KjSh:449} deals with unsuperstable f.o. theories both under elementary
embeddings and \cite{KjSh:455} deals with classes of abelian groups under pure
embeddings. 
\newline
If $\lambda = \lambda^{\aleph_0}$ there are universal groups of cardinality
$\lambda$ (among the reduced separable abelian p-groups): compact ones (see
\cite{Fu}, \cite{KjSh:455}).
\newline
If $2^{\aleph_0} < \lambda = \text{cf}(\lambda)$ and $\mu^+ < \lambda < 
\mu^{\aleph_0}$ then by Kojman Shelah \cite{KjSh:455} we have non-existence
results for \underbar{pure} embeddings.  Here we get the results for any
embedding, restricting ourselves to $(< \lambda)$-stable groups (see
Definition 3 below).  Without this restriction we shall deal with it in
\cite{Sh:552}. 
\newline 
If $\lambda < 2^{\aleph_0}$ on independence results for existence see
\cite{Sh:552}.  More results in these directions may be found in
\cite{Sh:457}, \cite{Sh:500}. 
\newpage

\head {\S1} \endhead
\bigskip

\demo{1.1 Context}  Fix $\bar n = \langle n_i:i < \omega \rangle$ a sequence
of natural numbers $> 1$.  We shall deal only with abelian groups, so we may
omit ``abelian".  \newline
For an (abelian) group $G$ define a prenorm $\| x\| = \text{ min}\{ 2^{-i}:x
\text{ divisible by } \dsize \prod_{j < i} n_j\}$ in $G$.  Let $d$ be the
induced semi-metric; i.e. distance is $d(x,y) = \| x-y \|$.  When in doubt use
$d_{\bar n},\| x \|_{\bar n}$. 
\enddemo
\bigskip

\demo{1.2 Notation}  Let ``group'' mean ``abelian group''. \newline
We concentrate on the classes 
\newline
(1)  ${\frak K}^{\text{rs}}_p = \{G:G \text{ is a } p\text{-group which is
reduced and separable}\}$ \newline
(where $G$ is separable if every pure subgroup of rank $1$ is a direct summand,
\newline
$G$ is a $p$-group means $(\forall x \in G)(\exists n < \omega)[n \ge 1 \and
p^n \, x = 0]$ and the norm here is $\| x\|_{<p:i < \omega >}$. 
\newline
(2) ${\frak K}^{\text{tfr}}_{\bar n} = \{G:G \text{ is a torsion free group
and }$, \newline
$\qquad \qquad \qquad d_{\bar n} \text{ is a metric (equivalently:}
\|-\|_{\bar n}  \text{ is a norm)}\}$. \newline

We say that $G$ is $\bar n$-reduced, if $G \in {\frak K}^{\text{tfr}}_{i_{\bar
n}}$; reduced means $\langle i:i < \omega \rangle$-reduced. \newline
(3)  Let ${\frak K}_\lambda = \{G \in {\frak K}:\| G \| = \lambda \}$.
\enddemo
\bigskip

\definition{1.3 Definition}  We say that $G$ is $(< \lambda)$-stable if:
\newline
$A \subseteq G,|A| < \lambda \Rightarrow \text{ closure}_G
(\langle A \rangle_G) = c \ell_G(\langle A \rangle_G) =
\{x:d(x,\langle A \rangle_G) = 0 \}$ \newline
has cardinality $< \lambda$
(where $\langle A \rangle_G =$ the subgroup of 
$G$ generated by $A$ and \newline
$d(x,A) = \text{ inf}\{d(x,y):y \in A\}$).
\enddefinition
\bigskip

\remark{1.4 Remark}  $\| - \|$ is a norm on $G$ if $\| x \| = 0 \Rightarrow
x = 0$.  For torsion free groups this means no non-zero homomorphic image of
$\left< \biggl\{ {\frac 1 {\dsize \prod_{j<i} n_j}}:i < \omega \biggr\}
\right>$ is embeddable.
\endremark
\bigskip

\proclaim{1.5 Claim}  If $G$ is a strongly $\lambda$-free group of cardinality
$\lambda$ \underbar{then} $G$ is $(< \lambda)$-stable. \newline 
If $(\forall \mu < \lambda)(\mu^{\aleph_0} < \lambda)$ \underbar{then}
every \footnote{i.e. $G \in {\frak K}^{\text{rs}}_p$, closure under
$d_{\langle p:i < \omega \rangle}$ \underbar{or} $G \in {\frak K}^{\text{tfr}}
_{\bar n}$, closure under $d_{\bar n}$, see 2.}
$G$ is $(< \lambda)$-stable.
\endproclaim
\bigskip

\remark{Remark}  The notion ``$G$ is strongly $\lambda$-free" is well known.
It means: \newline
$(\forall H)(\exists K)[H \subseteq G \and |H| < \lambda
\rightarrow H \subseteq K \subseteq G \and K$ is free and $G/K$ is
$\lambda$ free $]$.  See \cite{Ek1}.
\endremark
\bigskip

\demo{1.6 Discussion}  Now in the results of \cite{KjSh:455}, we can consider
not necessarily pure embeddings.  A difference from there is that here we add
$$A_\delta \subseteq \{ \alpha \in C_\delta:cf(\alpha) > \aleph_0 \}$$
(see \cite{KjSh:455}, Lemma 2, clause (iii).)
\enddemo
\bigskip

\definition{1.6 Definition}  Suppose that $\lambda$ is a regular uncountable
cardinal and that $G$ is a group of cardinality $\lambda$. \newline
1)  A sequence
$\bar G = \langle G_\alpha:
\alpha < \lambda \rangle$ is called a $\lambda$-representation of $G$ if and
only if for all $\alpha$:
\medskip
\roster
\item $G_\alpha \subseteq G_{\alpha + 1}$
\item $G_\alpha$ is of cardinality smaller than $\lambda$
\item if $\alpha$ is limit then $G_\alpha = 
\dsize \bigcup_{\beta < \alpha} G_\beta$
\item $G = \dsize \bigcup_{\alpha < \lambda} G_\alpha$.
\endroster
\medskip

\noindent
2) Suppose $\bar G = \langle G_\alpha:\alpha < \lambda \rangle$ is a given
representation of a group $G$.  Suppose $c \subseteq \lambda$ is a set of
ordinals, and the increasing enumeration of $c$ is $\langle \alpha_i: i <
i(*) \rangle$.  Let $g \in G$ be an element.  We define a way in which $g$
chooses a subset of $c$:

$$
\align
\text{Inv}_{\bar G}(g,c) = \{ \alpha_i \in c:&\text{ for some }
n \in [2,\omega) \text{ we have } \\
  &g \in (G_{\alpha_i + 1} + nG) \text{ but } g \notin (G_{\alpha_i} +
  nG)) \}.
\endalign
$$

We call $\text{Inv}_{\bar G}(g,c)$ the invariant of the element $g$ relative
to the $\lambda$-representation $\bar G$ and the set of indices $c$.

Worded otherwise, $\text{Inv}_{\bar G}(g,c)$, is the subset of those indices
$\alpha_i$ such that by increasing the group $G_{\alpha_i}$ to the larger
group $G_{\alpha_{i+1}}$, an $n$-congruent for $g$ is introduced for some
$n$.
\medskip

\noindent
3) Suppose that $\bar C = \langle C_\delta:\delta \in S \rangle$ is
a club guessing sequence; i.e. $C_\delta \subseteq \delta$ and for every
club $E$ of $\lambda$ for stationarily many $\delta \in S$, we have $C_\delta
\subseteq E$; we do not require that $C_\delta$ is a club of $\delta$
and that $\bar G$ is a $\lambda$-representation of a group $G$ of
cardinality $\lambda$.  Let
\medskip
\roster
\item "{(a)}"  $P_\delta(\bar G,\bar C) = \{ \text{Inv}_{\bar G}
(g,C_\delta):g \in G\}$
\item "{(b)}"  $P'_\delta(\bar G,\bar C) = \{ \text{Inv}_{\bar G}(g,C_\delta)
:x \in G$, moreover $x \in c \ell_G(G_\delta)\}$
\item "{(c)}"  $\text{INV}(G,\bar C) = 
[ \langle P'_\delta(\bar G,\bar C):\delta \in S \rangle]/id(\bar C)$.
\endroster
\medskip

\noindent The second item should read ``the equivalence class of the sequence
of $P_\delta$'s modulo the ideal $\text{id}(\bar C)$", where two sequences are
equivalent modulo an ideal if the set of coordinates in which the sequences
differ is in the ideal and 
$$\text{id}(\bar C) = \text{id}^a(\bar C) = \{ A \subseteq \lambda:\text{ for
every club } E \text{ of } \lambda \text{ for no } \delta \in A \cap E,
C_\delta \subseteq E\}.$$
\enddefinition
\bigskip

\demo{1.7 Fact}  For $\lambda,G,\bar G,\bar C,C$ as above:
\medskip
\roster
\item $\text{Inv}_{\bar G}(y,C)$ is a countable subset of $C$
\item $P_\delta(\bar G,\bar C)$ is a family of countable subsets of
$C_\delta$ of cardinality $\le |C|^{\aleph_0}$ and $\le \lambda$
\item  $P'_\delta(\bar G,\bar C)$ is a subset of $P_\delta(\bar G,\bar C)$,
if $G \in {\frak K}^{\text{rs}}_p$ then $P'_\delta(\bar G,\bar C) = P_\delta
(\bar G,\bar C)$.
\item  If in addition $G$ is $(< \lambda)$-stable then $P'_\delta(\bar G,
\bar C)$ has cardinality $< \lambda$. 
\item $\text{INV}(G,\bar C)$ really does not depend on the choice of the
representation $\bar G$.
\endroster
\enddemo
\bigskip

\demo{Proof}  Straight.
\enddemo
\bigskip

\proclaim{1.8 Lemma}  Let $\lambda > 2^{\aleph_0}$ be regular, $\bar C =
\langle C_\delta:\delta \in S \rangle$ be a club guessing sequence, $\bar A =
\langle A_\delta:\delta \in S \rangle$, $A_\delta \subseteq
C_\delta,\text{otp}(A_\delta) = \omega$ and $\alpha \in A_\delta \Rightarrow
\text{cf}(\alpha) > \aleph_0$.
\newline
\underbar{Then} there is a $(< \lambda)$-stable separable $p$-group $G$,
and $|G| = A$ such that
\medskip
\roster
\item "{$(*)$}"  if $H$ is $(< \lambda)$-stable separable $p$-group $H$,
$|H| = \lambda,G$ embeddable into $H,\bar H$ a representation of $H$ 
\underbar{then} for $id(\bar C)$-almost every $\delta \in S$ we have
\item "{$(**)$}"  $A_\delta \subseteq B$ for some $B \in 
P_\delta(\bar H,\bar C)$.
\endroster
\endproclaim
\bigskip

\remark{Remark}  The proof is closely related to \cite{Sh:e}, III, 7.15.
\endremark
\bigskip

\demo{Proof}  Let for $\delta \in S,\eta_\delta$ be an increasing
$\omega$-sequence of ordinals enumerating $A_\delta$.  Let 
\newline
$\Gamma_n = \biggl\{ \rho:\rho \in {}^n \omega,\rho \text{ strictly
increasing sequence of length } n,\rho(\ell) \in \omega \backslash
\{ 0 \} \biggr\}$ \newline
and let $\Gamma = \dsize \bigcup_{n < \omega} \Gamma_n$ and
$\Gamma_\omega = \{\rho:\rho \in {}^\omega \omega$ and $n < \omega
\Rightarrow \rho \restriction n \in \Gamma_n\}$. 
\newline 
Let $G^*$ be generated as a group by $x_{\eta,\rho}$ where for some $n <
\omega$ we have $\eta \in {}^n \lambda$ and $\rho \in \Gamma_n$ freely except:

$$
p^{g(\rho)+1} x_{\eta,\rho} = 0 \text{ where } g(\rho) =: \sum \{ \rho(m):
m < \ell g(\rho) \}.
$$

Let $G$ be generated by $x_{\eta,\rho}$ (where $\dsize \bigvee_n(\eta \in {}^n
\lambda \and \rho \in \Gamma_n))$, $y_{\delta,\rho,n}$ ($\delta \in S$, $\rho
\in \Gamma$ and $n < \omega$) freely except:

$$
\oplus_1 \qquad \quad  p^{g(\rho)+1}x_{\eta,\rho} \equiv 0
\text{ (for } \eta \in {}^n \lambda \text{ and } \rho \in \Gamma_n)
$$

$$
\oplus_2 \qquad \quad  y_{\delta,\rho,n}-p^{\rho(n)}
y_{\delta,\rho,n+1} = x_{\eta_\delta \restriction n,\rho \restriction n}
\text{ (for }\delta \in S,\rho \in \Gamma_\omega).
$$
\medskip

Let $G_\alpha$ be the subgroup of $G$ generated by

$$
Y_\alpha = \{x_{\eta,\rho}:(\exists n)(\eta \in {}^n \alpha \and
\rho \in \Gamma_n \} \cup \{ y_{\delta,\rho,n}:\delta \in S \cap \alpha,\rho
\in \Gamma_\omega,n < \omega\}.
$$

Clearly $G$ is from ${\frak K}^{\text{tr}}_p$, has cardinality $\lambda$ and
$\bar G = \langle G_\alpha:\alpha < \lambda \rangle$ is a representation of
$G$.  In fact $G_\alpha$ is the group generated by $Y_\alpha$ freely
except the equations of the form $\oplus_1,\oplus_2$ required above which
involve only its generators.
\medskip

Suppose $f:G \rightarrow H$ is an embedding and $H$ (not necessarily pure)
$\bar H$ are in $(*)$ of the Lemma.
\newline
Let

$$
\align
E =: \biggl\{ \alpha < \lambda:&f(G_\alpha) \subseteq H_\alpha,f^{-1}
\left( H_\alpha
\cap \text{ Rang}(f) \right) \subseteq G_\alpha, \\
  &\alpha \text{ a limit ordinal and for every } \beta < \alpha \\
  &\text{for some } \gamma, \beta < \gamma < \alpha \text{ and } \\
  &c \ell_G(G_\beta) \subseteq G_\gamma,c \ell_H(H_\beta) \subseteq H_\gamma,
\text{ for every } x_{\eta,\rho} \in G \\
  &\text{ we have } [x_{\eta,\rho} \in G_\alpha \Leftrightarrow \eta \in
{}^{\omega >}\alpha] \biggr\}.
\endalign
$$

\medskip
\noindent
Hence
\roster
\item "{$\oplus_3$}"  if $\alpha \in E, y \in H \backslash H_\alpha$ and
$\text{cf}(\alpha) > \aleph_0$ \underbar{then} $d(g,H_\alpha) > 0$
\item "{$\oplus_4$}"  if $\alpha \in E,\text{cf}(\alpha) > \aleph_0,
(\exists n)
[\eta \in {}^n \lambda \and \rho \in \Gamma_n] \text{ but } \eta \notin
{}^{\omega >} \alpha,x_{\eta,\rho} \in G$ \underbar{then}: \newline
$p^{g(\rho)}x_{\eta,\rho} \in G \backslash G_\alpha$ hence $p^{g(\rho)}
f(x_{\eta,\rho}) = f(p^{g(\rho)}x_{\eta,\rho}) \in H \backslash H_\alpha$
hence 
\newline
$d(p^{g(\rho)}x_{\eta,\rho},H_\alpha) > 0$.
\endroster
\medskip

\noindent
Let $\delta \in E$ be such that $C_\delta \subseteq E$.  Define a function
$h:\Gamma \rightarrow \omega$ such that:

$$
\qquad \quad (*)_0 \qquad \text{for } 
\rho \in {}^{n+1}\omega \text{ we have } 
d \left( f(p^{g(\rho)}x_{\eta_\delta \restriction
(n+1),\rho}),H_{\eta_\delta(n)} \right) \ge 2^{-h(\rho)}.
$$

\noindent
Note: $h$ is defined by $\oplus_4$. 
\newline
Next choose a  strictly increasing function  $\rho^* \in {}^\omega \omega$
such that: 

$$
\text{for every } n < \omega \text{ we have }
\rho^*(n) > h(\rho^* \restriction n)
$$

\noindent
(this is an overkill). Remembering $g(\eta) = \dsize \sum_{m < \ell g(\eta)} 
\rho^*(m)$ we clearly can prove by induction on $n$ that (use $\oplus_2$ and
the definition of $g$): 
\medskip
\roster
\item "{$(*)_1$}"  $y_{\delta,\rho^*,0} = \dsize \sum_{m<n} p^{g(\rho^*
\restriction m)}x_{\eta_\delta \restriction m,\rho^* \restriction m}
+ p^{g(\rho^* \restriction n)}y_{\delta,\rho^*,n}$.
\endroster
\medskip

\noindent
Note: $x_{\eta_\delta \restriction (n+1),\rho^* \restriction (n+1)} \in
G_{\eta_\delta (n+1)} \backslash G_{\eta_\delta(n)}$ and $\eta_\delta(\ell)
< \eta_\delta(n)$ for $\ell < n$ hence
\newline
$\dsize \sum_{m<n} p^{g(\rho^* \restriction m)}
x_{\eta_\delta \restriction m,\rho^* \restriction m} \in G_{\eta_\delta}(n)$.
By this and $(*)_1$ for each $n < \omega$ (use $n+2$ above)
\medskip
\roster
\item "{$(*)_2$}"  $y_{\delta,\rho^*,0} - p^{g(\rho^* \restriction (n+1)}
x_{\eta_\delta \restriction (n+1),\rho^* \restriction (n+1)} \in
G_{\eta_\delta(n)} + p^{g(\rho^* \restriction (n+2)}G$
\endroster
\medskip

\noindent
hence (as $f$ is an embedding of $G$ into $H$ mapping $G_{\eta_\delta(n)}$
into $H_{\eta_\delta(n)}$ as \newline
$\eta_\delta(n) \in A_\delta \subseteq C_\delta \subseteq E)$:
\medskip
\roster
\item "{$(*)_3$}"  $f(y_{\delta,\rho^*,0}) - 
p^{g(\rho^* \restriction (n+1)}f(x_{\eta_\delta \restriction (n+1),\rho^*
\restriction (n+1)}) \in H_{\eta_\delta(n)} + p^{g(\rho^* \restriction
(n+2)}H$.
\endroster
\medskip

By the choice of $\rho^*$ we have

$$
\rho^*(n+1) > h(\rho^* \restriction (n+1))
$$

\noindent
and hence by the choice of $h$ (i.e. by $(*)_0$) we conclude
\medskip
\roster
\item "{$(*)_4$}"  $p^{g(\rho^* \restriction (n+1)}
f(x_{\eta_\delta \restriction (n+1),\rho^* \restriction (n+1)}) \notin
H_{\eta_\delta(n)} + p^{g(\rho^* \restriction (n+2)}H$.
\endroster
\medskip

\noindent
By $(*)_3$ and $(*)_4$ we get
\medskip
\roster
\item "{$(*)_5$}"  $f(y_{\delta,\rho^*,0}) \notin
H_{\eta_\delta(n)} + p^{g(\rho^* \restriction (n+2))}H$.
\endroster
\medskip

\noindent Now use $(*)_1$ for $\delta,\rho^*,n+2$: note that $\dsize \sum_{m <
n+2}p^{g(\rho^* \restriction m)}x_{\eta_\delta \restriction m,\rho^*
\restriction m}$ belongs to $G_\gamma$ if $\dsize \bigcup_{m < n+2}
\text{Rang}(n_\delta \restriction m) \subseteq \gamma \in E$, but the maximal
ordinal in $\dsize \bigcup_{m < n+2} \text{Rang}(\eta_\delta \restriction m)$ is $\eta_\delta(n)$; hence
\medskip
\roster
\item "{$(*)_6$}"  $y_{\delta,\rho^*,0} \in H_{\text{min}(C_\delta
\backslash (\eta_\delta(n) + 1))} + p^{g(\rho^* \restriction (n+2)}H$.
\endroster
\medskip

\noindent
Apply $f$ on $(*)_6$ as $C_\delta \subseteq E$ we have
\medskip
\roster
\item "{$(*)_7$}"  $f(y_{\delta,\rho^*,0}) \in
H_{\text{min}(C_\delta \backslash (\eta_\delta(n)+1)}
+ p^{g(\rho^* \restriction (n+2))}H$.
\endroster
\medskip

\noindent
So by $(*)_5 + (*)_7$ we deduce $\eta_\delta(n) \in \text{Inv}_{\bar H}
\left( f(y_{\delta,\rho^*,0}),C_\delta \right)$ for each
$n$, so \newline
$A_\delta \subseteq \text{Inv}_{\bar H} \left(
f(y_{\delta,\rho^*,0}),C_\delta \right)$, but the later belongs to
$P_\delta(\bar H,C_\delta)$ as required. \hfill$\square_{1.8}$
\enddemo
\bigskip

\demo{1.9 Conclusion}  Assume $\lambda = cf(\lambda) > 2^{\aleph_0},
\mu^+ < \lambda < \mu^{\aleph_0}$. \underbar{Then} in the class

$$
\align
{\frak K}^{\text{rs},\text{st}}_{p,\lambda} =: \{ G:\ &G \text{ is an abelian
reduced, separable } p\text{-group, and} \\
  &G \text{ is } (< \lambda) \text{-stable of cardinality } \lambda \}
\endalign
$$

\noindent
there is no universal member under (usual) embedding.
\enddemo
\bigskip

\demo{Proof}  By \cite{Sh:420}, 1.8 we can find stationary $S \subseteq
\lambda$ and $\bar C = \langle C_\delta:\delta \in S \rangle,C_\delta$ a
subset of $\delta,\text{otp}(C_\delta) = \mu \times \omega$ (in fact,
$C_\delta$ closed in $\text{sup}(C_\delta)$, which is not necessarily
$\delta$) such that $\text{id}^a(\bar C)$ is a proper ideal on $\lambda$ (i.e.
$\lambda \notin \text{ id}^a(\bar C)$): remember

$$
\text{id}^a(C) = \{ A \subseteq \lambda:\text{for some club }E \text{ of }
\lambda \text{ for some } \delta \text{ is } C_\delta \subseteq E \}.
$$

\noindent
Clearly $\mu > \aleph_1$ (as we can replace $C_\delta$ by any
$C'_\delta \subseteq C_\delta$ of order type $\mu$) so without loss 
of generality if
$\alpha \in \text{ nacc}(C_\delta) =: \{ \beta \in C_\beta:\beta >
\text{ sup}(\beta \cap C_\beta)\}$ then $\text{cf}(\beta) > \aleph_0$.  
Suppose
$H \in {\frak K}^{\text{rs},\text{st}}_{p,\lambda}$ is universal in it, let
$\bar H = \langle H_i:i < \lambda \rangle$ be a representation of $H$.  Let

$$
P_\delta =: \{ \text{Inv}_{\bar H}(y,C_\delta):y \in H\}
$$

$$
P'_\delta =: \{A \subseteq C_\delta:\text{for some }B \text{ we have }
A \subseteq B \in P_\delta \}.
$$

So $P'_\delta$ is a family of $\le \lambda$ subsets of $C_\delta$, hence
there is $A_\delta \subseteq \text{ nacc}(C_\delta)$, unbounded in $C_\delta$
or order type $\omega$ such that $A_\delta \notin P'_\delta$.  Now
apply Lemma 1.8. \hfill$\square_{1.9}$ 
\enddemo
\bigskip

\proclaim{1.10 Claim}  In Lemma 1.8, instead of using $\Gamma_\omega$ use any
$\Gamma' \subseteq \Gamma_\omega$ such that:
\medskip
\roster
\item "{$\otimes_1$}"  $(\forall h)(h$ a function from
${}^{\omega >}\omega \text{ to } \omega \Rightarrow$ \newline
$\qquad(\exists \rho^* \in \Gamma')(\forall^\infty n)
(\rho^*(n) > h(\rho^* \restriction n))$
\endroster
\medskip

\noindent
is enough. \newline
So $|\Gamma| = {\frak d}$ is O.K.
\endproclaim
\bigskip

\demo{Proof}  Reflect.
\enddemo
\bigskip

\proclaim{1.11 Claim}  1) We can phrase Lemma 1.6 by an invariant, letting
\newline
$C_\delta = \{ \alpha_i:i < \text{ otp}(C_\delta)\}$ (increasing) by
defining:

$$
\align
P(x,C_\delta,\bar G) =: 
\biggl\{ A \subseteq C_\delta:&\text{otp}(A) = \omega,
A=\{\alpha_{i(n)}:n < \omega\},i(n) < i(n+1), \\
  &\text{every } \alpha \in A \text{ has cofinality } > \aleph_0, 
\text{ and there are} \\
  &T \subseteq {}^{\omega >}\omega, \text{ downward closed (by }
\triangleleft),<> \in T \text{ such that } \\
  &\eta \in T \Rightarrow (\exists^\infty n)[\eta \char 94 <n> \in T]
\text{ and } x_\rho \in G_{\alpha_{i(n)+1}} \text{ has} \\
  &\text{order } g(\rho) + 1,\text{ and for some } 
h:T \rightarrow \omega, \text{ for every} \\
  &\rho \in \text{ Lim}(T)(\subseteq {}^\omega \omega)
\text{ satisfying } (\forall n)[\rho(n) > h(\rho \restriction n)] \\
  &\text{there is } y_\rho \text{ such that } \\
  &y_\rho = \dsize \sum_{\ell < n} p^{g(\rho \restriction \ell)} 
x_{\rho \restriction \ell} \text{ mod } p^{g(\rho \restriction (n+1))}G
\text{ for each } n \biggr\}.
\endalign
$$
\medskip

\noindent 
2)  We can restrict ourselves to $\Gamma' \subseteq \Gamma_\omega$ as in
Claim 1.10.
\endproclaim
\bigskip

\proclaim{1.12 Claim}  1) We can deal similarly with the class
${\frak K}^{\text{tf}}_{\bar n}$, i.e. this class has no universal member in
$\lambda$ (under usual embedding if $2^{\aleph_0} < \mu,\mu^+ < \lambda =
\text{cf}(\lambda) < \mu^{\aleph_0}$.  In particular the class of reduced
torsion free
groups (using $\bar n = \langle i!:i < \omega \rangle$,etc). \newline
2)  Similarly, for
the class of metric spaces with embedding meaning one to one functions
preserving $\lim x_n = x$. \newline
3)  If in Lemma 1.8 we weaken $(**)$ to
\medskip
\roster
\item "{$(**)^-$}"  $A_\delta \cap B$ is infinite for some $B \in P_\delta
(\bar H,\bar C)$, then this:
{\roster
\itemitem{ (a) }  suffices in the proof of conclusion (a)
\itemitem{ (b) }  in order that it holds, the requirement in 1.10(1) can 
be weakened to:
\endroster}
\item "{$\otimes_2$}"  for every function $h$ from ${}^{\omega >}\omega$ to
$\omega$ we have \newline
$(\exists \rho^* \in \Gamma')(\exists^\infty n)[\rho^*(n) > h(\rho^*
\restriction n)]$.
\endroster
\medskip

\noindent
4)  If
\roster
\item "{$(a)$}"  $\lambda = \text{cf}(\lambda) > \min\{ \Gamma' \subseteq 
\Gamma_\omega:\Gamma' \text{ satisfies } \otimes_1\}$ or just
\item "{$\quad(a)^-$}"  $\lambda = \text{cf}(\lambda) > \min\{|\Gamma'|:\Gamma'
\subseteq \Gamma_\omega,\Gamma' \text{ satisfies } \otimes_2\}$ and
\item "{$(b)$}"  for some $\mu < \lambda$ we have $cf([\mu]^{\aleph_0},
\subseteq) > \lambda$ \newline
\endroster
\underbar{then} ${\frak K}^{\text{rs},\text{st}}_{p,\lambda}$ has no universal
member. 
\medskip

\noindent
(in the above definitions $(a)$ or $(a)^-$ we can replace $\Gamma_\omega$
by ${}^\omega \omega$, does not matter).
\endproclaim
\bigskip

\demo{Proof}  E.g. \newline
(3)(a)  Just note that there is $\langle B^\delta_i:i < \mu^{\aleph_0}
\rangle,B^\delta_i \subseteq C_\delta$ is unbounded in it and $i \ne j
\Rightarrow B^\delta_i \cap B^\delta_j$ is finite.  So every $B \in P_\delta$
(see the proof of Conclusion 1.9) we have $\{ i < \mu^{\aleph_0}:B \cap
B^\delta_i$ is infinite$\}$ has cardinality $\le 2^{\aleph_0}$. \newline
(3)  The point is that we can find an order $<^*$ on $\mu \times \omega$
such that $(\mu \times \omega,<^*)$ is a tree with $\omega$ levels the
$(n+1)$-th level is $(\mu \times n,\mu \times (n+1))$, level $0$ is
$\{ 0 \}$, level $1$ is $[1,\mu)$, and each node has $\mu$ immediate
successor, we can have
$g:\mu \times \omega \rightarrow [\mu \times \omega]^{\aleph_0}$ such that 
$w \in [\mu \times \omega]^{< \aleph_0} \Rightarrow \mu \times \omega
= \text{ otp}\{\alpha:g(\alpha) = w\}$. \newline
If ${\Cal P} \subseteq [\mu \times \omega]^{\aleph_0}$, $|{\Cal P}| <
\text{cf}([\mu]^{\aleph_0},\subseteq)$ without loss of generality
$$B \in {\Cal P} \and \alpha \in B \Rightarrow g(\alpha) \subseteq B.$$ 
Choose $A \in [\mu \times \omega]^{\aleph_0}$ not included in any $B \in A$,
let $A = \{\alpha_n:n < \omega\}$, choose $\gamma_n,\gamma_{n+1} > \gamma_n +
\mu,g(\gamma_n) = \{\alpha_\ell:\ell < n\}$. 
\enddemo
\newpage

\head {\S2 Large cofinality and singulars} \endhead
\bigskip

We want to prove that ${\frak K}^{\text{rs}}_{p,\lambda}$ has large
cofinality (parallel to \cite{KjSh:409} \S4).  We meet a problem here as we
have to deal with $y_{\delta,\rho,0}(\rho \in \Gamma_\omega)$ not just with
$y_{\delta,0}$. 
\newline
Remember: existence of universal is equivalent to cofinality being 1, i.e.
cofinality of $(({\frak K}^{\text{rs}}_p)_\lambda$, embedability) and, what is
closely related, prove non-existence of universal in singular cardinals.
\bigskip

\definition{2.1 Definition}  For an ideal $J$ on $\lambda$ and $\chi \ge
\lambda$ let

$$
\align
\chi^{\langle J \rangle} = \text{ Min} \{| {\Cal P}|:&{\Cal P} \subseteq
[\chi]^\lambda,\text{ and for every }f \in {}^\lambda \chi \\
 &\text{there is }B \in {\Cal P}, \text{ such that } \\
 &\{ i < \lambda:f(i) \in B \} \ne \emptyset \text{ mod } J \}
\endalign
$$

$$
\align
\chi^{[\lambda]} =
\text{ Min}\{ |{\Cal P}|:&{\Cal P} \subseteq [\chi]^\lambda,
\text{ and every } A \in [\chi]^\lambda \\
  &\text{is included in the union of } < \lambda \text{ members of }
{\Cal P} \}.
\endalign
$$
\enddefinition
\bigskip

\definition{2.2 Definition}  For $\lambda = \text{ cf}(\lambda) > \aleph_0$, 
and a club guessing sequence $\bar C = \langle C_\delta:\delta \in S \rangle$
we define $J = J^\ell_*[\bar C]$ (for $\ell = 1,2$).  It is the
following family of subsets of $\lambda \times ({}^\omega \omega)$: \newline

$$
\align
\biggl\{ A \subseteq \lambda \times ({}^\omega \omega):&\text{ for some }
B \in \text{id}^a[\bar C] \text{ for every} \\ 
  &\delta \in S \backslash B \text{ there is } h:{}^{\omega >} \omega
\rightarrow \omega \text{ such that for no }\rho \\
  &\text{do we have: } (\delta,\rho) \in A \text{ and:} \\
  &\ell = 0 \Rightarrow \text{ for infinitely many } n < \omega,
\rho(n) > h(\rho \restriction n) \\
  &\ell = 1 \Rightarrow \text{ for every large enough } n, \\
  &\rho(n) > h(\rho \restriction n) \biggr\}.
\endalign
$$

\enddefinition  
\bigskip

\proclaim{2.3 Claim}  Assume that $\lambda = \text{cf}(\lambda) >
2^{\aleph_0}$, $\bar C =\langle C_\delta:\delta\in S\rangle$ is a club
guessing sequence, $S \subseteq \lambda$, $\mu^+<\lambda\le\chi_0\le\chi_1<
\mu^{\aleph_0}$, $H_j\in {\frak K}^{\text{rs},\text{st}}_{p,\le \chi_0}$ for
$j < \chi_1$, and $\chi^{[J^0_*[\bar C]]} < \mu^{\aleph_0}$.  \underbar{Then}
there is $G \in {\frak K}^{\text{rs},\text{st}}_{p,\lambda}$ which is not not
embeddable in $H_\gamma$ for every $\gamma < \chi_1$.
\endproclaim
\bigskip

\demo{Proof}  Without loss of generality $|H_\gamma| = \chi_0$ and the set
of elements of $H_\gamma$ is $\chi_0$, let $\chi^*_0 =
\chi_0^{[J^1_*[\bar C]]}$ and let
${\Cal P} = \{A_\varepsilon:\varepsilon < \chi^*_0 \}$ exemplifies
$\chi^{[J^1_*[\bar C]]} = \chi^*_0$.
Let $H_{\gamma,\varepsilon}$ be a pure subgroup of
$H_\gamma$ of cardinality $\lambda$ including $A_\varepsilon$ and
$\bar H_{\gamma,\varepsilon} = \langle H_{\gamma,\varepsilon,i}:i < \lambda
\rangle$ be a representation of $H_{\gamma,\varepsilon}$.  Let
$\bar C = \langle C_\delta:\delta \in S \rangle$ be as in the proof of 
Lemma 1.8, and
choose $A_\delta \subseteq \text{ nacc}(C_\delta)$ unbounded of order type
$\omega$ such that for no $\gamma < \chi$, and $\varepsilon < 
\chi^*_0$ is there a set $B$ such that $B \in P_\delta(\bar H
_{\gamma,\varepsilon},\bar C)$ and $A \cap B$ infinite
(possible by cardinality
considerations).  Let $G$ be as in Lemma 1.8 for $\lambda,\bar C,
\langle A_\delta:\delta \in S \rangle$.  If $G$ is embeddable into
say $H_\gamma$, let $f$ be such embedding, so for some
$\zeta < \chi^*_0$ for $\zeta < \zeta(*)$ we have

$$
\{ (\delta,\rho):\delta \in S,\rho \in {}^\omega \omega \} \notin
J^0_*[\bar C]\}.
$$

Now repeat the proof of Lemma 1.8. \hfill$\square_{2.3}$
\enddemo
\bigskip

\demo{2.4 Conclusion}  Assume $\lambda=\text{cf}(\lambda)>2^{\aleph_0},\mu^+
< \lambda < \chi_0 \le \chi_1 < \mu^{\aleph_0}$ and
$\chi^{[J^0_*[\bar C]]}_0 < \mu^{\aleph_0}$ for some guessing club system
$\bar C = \langle C_\delta:\delta \in S \rangle$, $S \subseteq \lambda$
stationary.

Then for every cardinal $\lambda'$ (possibly singular) $\lambda \le \lambda'
\le \chi_0$, the class ${\frak K}^{\text{rs},\text{st}}_{p,\lambda}$ has no
universal member (under usual embeddings).
\enddemo
\bigskip

\demo{2.5 Question}  Can we combine $1.10(4) + 1.2.4$?
\enddemo
\bigskip

\centerline {$* \qquad * \qquad *$}
\bigskip

\remark{Remark} The above can also be interpreted as giving negative
results on the existence of models of cardinality $\chi_0 > \lambda$
universal for ${\frak K}^{\text{rs},\text{st}}_{p,\lambda}$, etc. For more
results on non existence of universals  see \cite{Sh:552}. 
\endremark
\vfill
\newpage

REFERENCES
\bigskip

\bibliographystyle{literal-unsrt}
\bibliography{lista,listb,listx}

\shlhetal
\enddocument
\bye